\documentclass[12pt]{article}
\usepackage{graphicx,amsmath,amssymb,amsthm} 

\title{Cyclic group representations for relation algebras $57_{65}$ and $63_{65}$}
\author{Jeremy F.~Alm\\ Southern Illinois University}
\date{April 2026}

\begin{document}

\maketitle

\begin{abstract}
   We exhibit finite cyclic group representations for relation algebras $57_{65}$ and $63_{65}$. As a consequence, of the ten symmetric integral RAs on four atoms having at least one flexible atom, all are now known to have a representation over a finite cyclic group except for $33_{65}$, which is not even known to be finitely representable. 
\end{abstract}

\section{Introduction}
In a finite integral relation algebra (FIRA), an atom $a$ is called \emph{flexible} if there are no forbidden diversity cycles involving $a$. Every FIRA with a flexible atom is representable over a countably infinite set. The Flexible Atom Conjecture (FAC) asserts that every such RA is actually representable over a finite set. 

While a proof of the FAC seems out of reach at the moment, there is a years-long project underway to eliminate potential counterexamples. Many finite representations have been found, any many are finite group representations.  Among symmetric FIRAs with at most four atoms and at least one flexible atom, all are known to have finite representations, save for the stubborn $33_{65}$. (See \cite{Maddux} for the numbering system.) Most are known to have finite group representations, and most also have cyclic group representations.  Some recent additions to the list of known group representations are $31_{37}$, representable over $\mathbb{Z}/ 33791\mathbb{Z}$, and $32_{65}$, representable over $\mathbb{Z}/ 751181\mathbb{Z}$ \cite{Monk};  $59_{65}$, representable over $\mathbb{Z}/ 113\mathbb{Z}$ \cite{twosmall}; $33_{37}$ representable over $\mathbb{Z}/ 29\mathbb{Z}$, and $35_{37}$, representable over $\mathbb{Z}/ 3221\mathbb{Z}$ \cite{directed}; and $1896_{3013}$, representable over $\mathbb{Z}/1531\mathbb{Z}$ \cite{Jas}. 

All over the cyclic representations in the previous paragraph were found using the method initiated in \cite{Comer} and made substantially more powerful via an algorithmic speed-up in \cite{fast}. 

In this paper, we give cyclic representations for $57_{65}$ and $63_{65}$. For the former, it is the first known finite group representation; for the latter, it is the first known \emph{cyclic} representation (to the best of our knowledge). 

It is noteworthy that for small RAs, cyclic representations seem abundant. For example, all Ramsey relation algebras that are known to be representable have cyclic representations \cite{401,manske,Kowalski}. See also \cite{cycspec}, summarized in Table \ref{tab:summary}, where the spectrum for cyclic representations was determined for all symmetric integral 3-atom algebras. (The ordinary spectra were determined in \cite{AndMad}.)

\begin{table}[h]
    \centering
    \begin{tabular}{r|c|c}
         & Spec & Cyclic Spec\\
     \hline
    $1_7$ & $\{4\}$ & $\{4\}$\\
    $2_7$ & $\{n \geq 6\}$ & $\{2k : k \geq 3\}$\\
    $3_7$ & $\{2k : k \geq 3\}$ & $\{2k : k \geq 3\}$\\
    $4_7$ & $\{ n \geq 9\}$ & $\{n \geq 9\}\setminus\{p, 2p : p \text{ prime} \}$\\
    $5_7$ & $\{5\}$ & $\{5\}$\\
    $6_7$ & $\{n \geq 8\}$ & $\{8\} \cup \{ n \geq 11 \}$\\
    $7_7$ & $\{ n \geq 9\}$ & $\{n \geq 12\}$\\
    
    \end{tabular}
    \caption{Summary of Results from \cite{AndMad} and \cite{cycspec}}
    \label{tab:summary}
\end{table}

\section{Cyclic Representations}
Relation algebra $63_{65}$ has two forbidden cycles, $bbb$ and $ccc$. We give a representation over $\mathbb{Z}/ 29\mathbb{Z}$. 
Let

\begin{align*}
    A &= \{3, 7, 8, 9, 11, 13, 16, 18, 20, 21, 22, 26\}\\
    B &= \{1, 4, 10, 12, 17, 19, 25, 28\}\\
    C &= \{2, 5, 6, 14, 15, 23, 24, 27\}
\end{align*}

This coloring, which also serves as proof that the Ramsey number $R(4,3,3)$ is greater than 29, was noticed by the author in \cite{Rad}, where it is attributed to \cite{Kalb}. (Note that for Ramsey theory purposes, only the forbidden triangles are relevant.)

Relation algebra $57_{65}$ has forbidden cycles $ccc$ and $cbb$.  We give a representation over $\mathbb{Z}/ 46\mathbb{Z}$. 

Let

\begin{align*}
    A &= \{1, 2, 9, 10, 12, 13, 15, 18, 20, 21, 22, 23, 24,\\ 
    &\phantom{=} \ \ 25, 26, 28, 31, 33, 34, 36, 37, 44, 45\}\\
    B &= \{6, 7, 8, 14, 16, 30, 32, 38, 39, 40\}\\
    C &= \{3, 4, 5, 11, 17, 19, 27, 29, 35, 41, 42, 43\}
\end{align*}

This representation was found with a SAT solver, but can be verified to be correct without a solver, and has been checked independently by Roger Maddux.  

\section{Conclusion}

Of the symmetric integral RAs on four atoms, ten have at least one flexible atom.  All but $33_{65}$ were known to be finitely representable, and now all but $33_{65}$ are known to be representable over a finite cyclic group.  

The forbidden cycles for $33_{65}$ are $ccc$, $bcc$, and $cbb$. The symmetry in the forbidden 2-cycles, but asymmetry in the forbidden 1-cycles, presents a challenge.  

The author has checked that $33_{65}$ is not representable over $\mathbb{Z}/ n\mathbb{Z}$ for any $n\leq 100$, nor the symmetric group $S_n$ for $n\leq 5$. 

\section*{Acknowledgements}
The author wishes to thank Roger Maddux, for 20 years of conversation about the flexible atom conjecture, and his current MS student Eli Atkins, for his energy and fresh ideas that have given this university administrator a renewed sense of vigor. 

\section{Declarations}
The author has no competing interests to declare that are relevant to the content of this article.

All relevant data are included in this article. 

\bibliographystyle{plain}
\bibliography{refs}

\end{document}